\documentclass{amsart}

\usepackage[mathcal]{euscript}
\usepackage{amssymb, amsfonts, xypic}

\newtheorem{thm}{Theorem}[section]
\newtheorem{cor}[thm]{Corollary}
\newtheorem{prop}[thm]{Proposition}
\newtheorem{lem}[thm]{Lemma}

\newtheorem{ques}[thm]{Question}

\theoremstyle{definition}
\newtheorem{defn}[thm]{Definition}
\newtheorem{ex}[thm]{Example}

\theoremstyle{remark}
\newtheorem{rem}[thm]{Remark}

\numberwithin{equation}{section}

\newcommand{\Map}{\textup{Map}}
\newcommand{\Fin}{\textup{Fin}}

\newcommand{\Ho}{\textup{Ho}}
\newcommand{\Hom}{\textup{Hom}}
\newcommand{\sk}{\textup{sk}}

\newcommand{\cosk}{\textup{cosk}}

\newcommand{\pro}{\textup{pro-}}
\newcommand{\pros}{\textup{pro}}

\newcommand{\map}{\rightarrow}
\renewcommand{\smash}{\wedge}
\newcommand{\bd}{\partial \Delta}

\newcommand{\C}{\mathcal{C}}
\newcommand{\F}{\mathcal{F}}

\newcommand{\Z}{\mathbb{Z}}
\newcommand{\N}{\mathbb{N}}

\newcommand{\A}{\mathcal{A}}
\newcommand{\bA}{\mathbb{A}}

\renewcommand{\H}{\mathcal{H}}

\newcommand{\ccomp}[1]{#1^{\wedge}_{R\mbox{-}c}}
\newcommand{\hcomp}[1]{#1^{\wedge}_{R\mbox{-}h}}
\newcommand{\tl}{\tilde}

\newcommand{\colim}{\mathop{\textup{colim}}}

\newcommand{\Ext}{\mathop{\textup{Ext}}}

\newcommand{\sSet}{\textit{sSet}}

\newdir{ >}{{}*!/-11pt/@{>}}
\newcommand{\cofib}{\ar@{ >->}}
\newcommand{\fib}{\ar@{->>}}

\newcommand{\op}{\text{op}}

\newcommand{\dfn}{\textbf} 
\newcommand{\mdfn}[1]{\dfn{\mathversion{bold}#1}} 

 \newcommand{\lab}[1]{\label{#1}}

\title{Completions of Pro-Spaces}
\author{Daniel C.\,Isaksen}
\address{
Department of Mathematics \\
Wayne State University \\
Detroit, MI 48202, USA   }

\thanks{This work was partially supported by a National Science
Foundation Postdoctoral Research Fellowship.  The author acknowledges
useful conversations with Bill Dwyer and Daniel Biss.}

\keywords{pro-space, singular cohomology, completion, Bousfield-Kan tower}

\subjclass{55P60, 55N10 (Primary); 18G55, 55U35 (Secondary)}
 
\begin{document}

\begin{abstract}
For every ring $R$, 
we present a pair of model structures on the category of pro-spaces.
In the first,
the weak equivalences are detected by cohomology with coefficients
in $R$.
In the second,
the weak equivalences are detected by cohomology with coefficients
in all $R$-modules (or equivalently by pro-homology with coefficients in $R$).
In the second model structure, fibrant replacement is essentially just
the Bousfield-Kan $R$-tower.
When $R = \Z/p$, the first homotopy category is equivalent to a 
homotopy theory defined by Morel but has some convenient categorical 
advantages.
\end{abstract}

\maketitle

\section{Introduction}
\label{sctn:intro}

The notion of $R$-completion has been a valuable tool to homotopy theorists.
The basic idea is to start with a space $X$ and then construct another
space $X^\wedge_R$ whose homotopy type is entirely determined by 
the singular 
cohomology $H^*(X; R)$ of $X$ with coefficients in $R$.  In other words,
$X^\wedge_R$ remembers the $R$-cohomology of $X$ but forgets all other
information.

Bousfield and Kan \cite{BK} constructed $R$-completions for a large class
of spaces that they called $R$-good spaces.  The basic construction 
goes as follows.  Start with a space $X$.  Then define a cosimplicial
space $\tilde{R}^* X$.  This cosimplicial space gives rise to a tower
\[
\cdots \map R_2 X \map R_1 X \map R_0 X
\]
of fibrations.  Finally, $X^\wedge_R$ is the homotopy limit of this tower.

Unfortunately, this process only works for the $R$-good spaces.  In fact,
the {\em tower} described above is correct for all spaces, but the
homotopy limit causes problems when $X$ is not $R$-good.  

Thus, one approach to generalizing the construction of Bousfield and Kan 
to arbitrary spaces is to consider $X^\wedge_R$ not as a single space but
rather as the whole tower.  In fact, it is best to think of $X^\wedge_R$
as a pro-space \cite{D}.
This paper is concerned with the homotopical foundations for
pro-spaces suitable for this viewpoint on $R$-completion.

When $R = \Z/p$, Morel \cite{Mo} constructed a homotopy theory
of simplicial pro-finite sets that is suitable
for studying $\Z/p$-completions of spaces.
Unfortunately, there are a few problems with the approach in \cite{Mo}.
Namely, it is {\em not} true that the category of simplicial pro-finite
sets is equal to the category of pro-simplicial finite sets.  In fact,
the former is only a retract of the latter.  This means that we must
be very careful with our intuitive ideas about simplicial pro-finite sets.

We prove the following theorem in Section \ref{sctn:cohlgy-mc}.

\begin{thm}
\label{thm:R-cohlgy-mc}
Let $R$ be any ring. 
There is a model structure on the category of pro-simplicial sets
in which 
the weak equivalences are maps $f: X \map Y$ such that
$H^*(Y;R) \map H^*(X; R)$ is an isomorphism.
\end{thm}

One of our main results is that the homotopy theory of 
pro-simplicial sets from Theorem \ref{thm:R-cohlgy-mc}
is equivalent to Morel's homotopy theory
of simplicial pro-finite sets when $R = \Z/p$.  
The advantage of our approach is that one does not have to worry about
the unintuitive nature of simplicial pro-finite sets.

The proofs of \cite{Mo} rely in an essential way
on notions of finiteness and take great advantage of the fact that
$\Z/p$ is a finite ring.  In fact, finiteness is not really such an important
ingredient, as demonstrated by the proof of Theorem \ref{thm:R-cohlgy-mc},
which works for any ring $R$.

In some contexts \cite{AM} \cite{I4}, one wants to take $R$-completions of
pro-spaces.  Although \cite{Mo} handles $\Z/p$-completions of spaces perfectly
well, it is not quite right for $R$-completions of arbitrary pro-spaces.
Here is the underlying reason.  If $X \map Y$ is a map of pro-spaces
such that $H^*(Y;R) \map H^*(X;R)$ is an isomorphism, then 
$H^*(Y;M) \map H^*(X;M)$ is {\em not} necessarily an isomorphism for
all $R$-modules $M$ (see Example \ref{ex}).  
These observations motivate the following theorem, which is also proved
in Section \ref{sctn:cohlgy-mc}.

\begin{thm}
\label{thm:R-hlgy-mc}
Let $R$ be any ring.
There is a model structure on the category of pro-simplicial sets
in which 
the weak equivalences are maps $f: X \map Y$ such that
$H^*(Y;M) \map H^*(X; M)$ is an isomorphism for all $R$-modules $M$.
\end{thm}

We will show later that the weak equivalences in the model structure
of Theorem \ref{thm:R-hlgy-mc} can also be described as the maps
$f: X \map Y$ such that $H_n(X;R) \map H_n(Y;R)$ is an isomorphism
of pro-groups for each $n \geq 0$.

For $R$-completions of pro-spaces, the
homotopy theory constructed in Theorem \ref{thm:R-hlgy-mc}
is better than the homotopy theory of Theorem \ref{thm:R-cohlgy-mc}.  
It retains more information, remembering not just
the cohomology with coefficients in $R$ but the cohomology
with coefficients in all $R$-modules.
Moreover, the homotopy theory of Theorem \ref{thm:R-hlgy-mc} has a close
link with the Bousfield-Kan $R$-tower \cite{BK}.
In particular, 
the Bousfield-Kan $R$-tower of a pointed connected space $X$ is
basically the same thing as 
the fibrant replacement of $X$ in the $R$-homological model structure.
We explore this link in Section \ref{sctn:R-comp}.

The main tool in the proof of Theorem \ref{thm:R-cohlgy-mc} is a 
general localization result from \cite{CI}.  The proof of 
Theorem \ref{thm:R-hlgy-mc} is similar, but we have to be careful about
the set-theoretical complications introduced by allowing cohomology
with coefficients in arbitrarily large $R$-modules.  

In both model structures, we give an explicit description of fibrant
objects.  This allows for computations, 
as demonstrated in Section \ref{sctn:free-gp},
where we compute the $R$-completion of the classifying space of a 
finitely-generated free group.
In order to describe the fibrant objects, we must introduce a notion
of nilpotence for spaces that is related to but distinct from the
usual notion of nilpotence.  See Section \ref{sctn:nilp} for more details.

At the end of the paper, we list several questions concerning this 
subject that remain unanswered.  We hope that this will encourage future
work on this topic.

\subsection{Background}

We assume familiarity with model structures.  The original reference
is \cite{Q}, but \cite{Ho} and \cite{Hi} are the modern thorough references.
Also, \cite{DS} is a good introduction to the subject.
We warn the reader about one convention concerning the model structure
axioms.  We will always start with a model category $\C$ in which 
factorizations are functorial.  However, when we produce model categories
on the pro-category $\pro \C$, we will not necessarily obtain functorial
factorizations.  Recent work of Chorny \cite{C} suggests that functoriality
is obtainable, but we will ignore the question here.

We work exclusively with simplicial sets, rather than topological spaces.
From now on, the word {\em space} always means {\em simplicial set}.
It is possible to obtain all of the results of this paper for topological
spaces instead of simplicial sets, but simplicial sets are somewhat 
easier to work with.  

We also assume a certain amount of familiarity with pro-categories,
although we give a brief review of the most important points in
Section \ref{sctn:pro-prelim}.  Although there are many established
and thorough references for pro-categories such as \cite{AM}, \cite{SGA4},
or \cite{EH}, the reader is encouraged to look also at \cite{I1}, \cite{I2},
and \cite{I3} for aspects of pro-categories that are particularly
relevant to this paper.

\subsection{Organization}

We begin in Section \ref{sctn:pro-prelim} with a brief review of
pro-categories, touching only on the issues that are most important
for present purposes.  
We also recall the strict model structure \cite{EH} \cite{I3}
on the category of pro-spaces and
a general localization result from \cite{CI} that will allow us
to produce new model structures for pro-spaces.
In the following section, we study the fibrant objects in these
localized model structures.  For this purpose, we need a variation on
the standard notion of nilpotence for spaces.

Section \ref{sctn:hom-alga} contains some relatively straightforward
material on homological algebra in pro-abelian categories; this is basically
just transferring well-known results to a new setting.

The main part of the paper begins in Section \ref{sctn:wk-eq} with
the study of pro-maps that induce isomorphisms in various kinds of
singular cohomology.  Section \ref{sctn:cohlgy-mc} describes the 
model structures that have these cohomology isomorphisms 
as weak equivalences.

The rest of the paper is dedicated to applications and connections to 
other established theories.  We begin in Section \ref{sctn:R-comp}
with the link between our constructions and the Bousfield-Kan
notion of $R$-completion.  Then in Section \ref{sctn:Morel}
we compare our constructions to Morel's theory of $\Z/p$-completion.
We warn the reader that Section \ref{sctn:Morel} is a bit tricky
because it involves a comparison of several categories that 
feel very similar but are definitely distinct.
In Section \ref{sctn:free-gp}, we 
compute the $\Z/p$-completion of the 
classifying space of a finitely-generated free group.
Finally, we list some open questions in
Section \ref{sctn:?}.

\section{Pro-Categories and Homotopy Theories for Pro-Spaces}
\label{sctn:pro-prelim}

\subsection{Pro-Categories}
\label{subsctn:pro-ctgy}

We begin with a brief overview of pro-categories and the homotopy
theory of pro-spaces.
Standard references on pro-categories include
\cite{SGA4}, \cite{AM}, and \cite{EH}. 
See also \cite{I2} and \cite{I3} for details specifically relevant
to the homotopy theory of pro-categories.

\begin{defn} 
\lab{defn:pro}
For a category $\C$, the category \mdfn{$\pro \C$} has objects all 
cofiltering
diagrams in $\C$, and 
$$\Hom_{\pro \C}(X,Y) = \lim_s \colim_t \Hom_{\C}
     (X_t, Y_s).$$
Composition is defined in the natural way.
\end{defn}

The word {\em pro-object} refers to objects of pro-categories.
A \dfn{constant} pro-object is one indexed
by the category with one object and one (identity) map.

A \dfn{level representation} of a map
$f:X \map Y$ is:
a cofiltered index category $I$;
cofiltered diagrams $\tl{X}$ and $\tl{Y}$ indexed by $I$ 
that are pro-isomorphic to $X$ and $Y$ respectively;
and a natural transformation $\tl{f}: \tl{X} \map \tl{Y}$
representing a pro-map that is isomorphic to $f$.
Every map has a level representation 
\cite[App.~3.2]{AM} \cite{Me}.

A pro-object $X$ satisfies a certain property \dfn{levelwise} if
each $X_{s}$ satisfies that property, and $X$ satisfies this property
\dfn{essentially levelwise} if it is isomorphic to another pro-object
satisfying this property levelwise.
Similarly, a level representation $X \map Y$ 
satisfies a certain property \dfn{levelwise}
if each $X_{s} \map Y_{s}$ has this property.
A map of pro-objects satisfies this property \dfn{essentially levelwise}
if it has a
level representation satisfying this property levelwise.

Let $c: \C \map \pro \C$ be the functor taking an object $X$ to the 
constant pro-object with value $X$.
Note that this functor makes $\C$ a full subcategory of
$\pro \C$.  The limit functor $\lim: \pro \C \map \C$
is the right adjoint of $c$.  
To avoid confusion, we write $\lim^{\pros}$ for limits computed 
within the category $\pro \C$.

We recall the construction of cofiltered limits in $\pro \C$ 
(see, for example, \cite[\S~4]{I1}).  The specific details of 
this construction will be used in several places later.
Start with a functor $X: A \map \pro \C: a \mapsto X^a$, where $A$ is a 
cofiltered index category.  
The index category $I$ for $\lim^{\pros} X$ consists of all pairs
$(a, s)$ such that $a$ belongs to $A$ and $s$ belongs to the indexing
category of $X^a$.  A morphism $(a,s) \map (b,t)$ consists of
a morphism $a \map b$ in $A$ together with a map $X^a_s \map X^b_t$ in
$\C$ that represents the pro-map $X^a \map X^b$.  Finally, 
$\lim^{\pros} X$ is defined to be the functor $I \map \C$ that
takes $(a,s)$ to $X^a_s$. 

\subsection{Strict Homotopy Theory of Pro-Spaces and Its Localizations}
\label{subsctn:pro-homotopy}

We now review from \cite{I3} the strict homotopy theory of 
pro-spaces.
The strict model structure was originally defined in \cite{EH}.  The 
\dfn{strict weak equivalences} ({\em resp.}, \dfn{cofibrations})
are the essentially
levelwise weak equivalences ({\em resp.}, cofibrations),
and the \dfn{strict fibrations} are defined by the right lifting property.
In fact, a more explicit description of the fibrations in
terms of matching maps is possible \cite[\S~4]{I3}.
The strict model structure is proper and simplicial.

Recall that the $n$th singular cohomology group $H^n(X; M)$ of a pro-space
$X$ with coefficients in an abelian group $M$ is defined to be
$\colim_s H^n(X_s; M)$
\cite[2.2]{AM} \cite{S}.
In fact, there is an isomorphism between $H^n(X;M)$ and the set
$[X, cK(M, n) ]_{\pros}$ of weak homotopy classes of maps 
of pro-spaces.
The pro-space $cK(M, n)$ is the constant pro-space 
with value an Eilen\-berg-Mac Lane space.

We recall the following localization result for pro-spaces.
The full proof (in greater generality),
which owes much to \cite[Ch.~5]{Hi}, appears in \cite{CI}.

\begin{thm}
\label{thm:local}
Let $K$ be a set of fibrant spaces.  There exists a 
left proper simplicial model structure
on the category of pro-spaces such that the cofibrations are
the essentially levelwise cofibrations and such that a map
$f:X \map Y$ is a weak equivalence if and only if
\[
\Map_{\pros}(Y, cA) = \colim_s \Map(Y_s, A) \map 
\colim_t \Map(X_t, A) = \Map_{\pros} (X, cA)
\]
is a weak equivalence for every object $A$ in $K$.
\end{thm}

The weak equivalences in the above theorem are called
\mdfn{$K$-colocal weak equivalences}.

\section{Nilpotent Spaces}
\label{sctn:nilp}

We will need to understand the 
fibrant objects
in the localized model categories of Theorem \ref{thm:local}.

\begin{defn}
\label{defn:K-nilp}
Let $K$ be any collection of fibrant spaces.  The class of 
\mdfn{$K$-nilpotent spaces}
is the smallest class of fibrant spaces such that:
\begin{enumerate}
\item the space $*$ is $K$-nilpotent;
\item the $K$-nilpotent spaces are closed under weak equivalences
between fibrant spaces;
\item and
if $X$ is $K$-nilpotent, $A$ belongs to $K$, and 
$X \map A^{\bd^k}$ is any map, then the 
fiber product $X \times_{A^{\bd^k}} A^{\Delta^k}$
is also $K$-nilpotent.
\end{enumerate}
\end{defn}

Because the fiber product in (3) above is actually a homotopy fiber
product, it is only the weak homotopy types of the spaces in $K$ that matter.
Thus, the $K$-nilpotent spaces are 
more properly a collection of weak homotopy types rather
than a collection of actual spaces.
The consequence is that we are allowed to choose any (fibrant) models
for the weak homotopy types in $K$ that are most convenient for our
purposes.

\begin{lem}
\label{lem:K-nilp}
A space $X$ is $K$-nilpotent if and only if it is fibrant
and weakly equivalent
to a space that can be built from $*$ by finitely many pullbacks of 
type (3) in Definition \ref{defn:K-nilp}.
\end{lem}

In other words, when constructing a $K$-nilpotent space, it is not necessary
to use any weak equivalences until the very last step.

\begin{proof}
Let $C_n$ be the class of $K$-nilpotent spaces 
that can be built from $*$ with fewer than $n+1$ pullbacks (and possibly
also weak equivalences), and let $D_n$ be the class of fibrant spaces
that are weakly equivalent to a space 
that can be built from $*$ with fewer than $n+1$ pullbacks (but without
any weak equivalences).
By definition, $D_n$ is contained in $C_n$.  
We will show by induction that $C_n$ and $D_n$ are equal.

The classes $C_0$ and $D_0$ consist of the fibrant contractible spaces,
so they are equal.  Now suppose that $C_{n-1}$ and $D_{n-1}$ are 
equal.  Let $X$ belong to $C_n$.  Then $X$ is weakly equivalent to
a space $X' \times_{A^{\bd^k}} A^{\Delta^k}$, where $X'$ belongs to
$C_{n-1}$.  By the induction assumption, $X'$ also belongs to $D_{n-1}$,
so $X'$ is weakly equivalent to a space $X''$ that can be built from $*$
by fewer than $n$ pullbacks.  Now $X' \times_{A^{\bd^k}} A^{\Delta^k}$
is weakly equivalent to $X'' \times_{A^{\bd^k}} A^{\Delta^k}$ because
the fiber products are actually homotopy fiber products.  
Thus
$X$ is weakly equivalent to $X'' \times_{A^{\bd^k}} A^{\Delta^k}$,
which is a space that can be built from $*$ using fewer than $n+1$
pullbacks.  Therefore, $X$ belongs to $D_n$.
\end{proof}

The next theorem demonstrates the
relevance of $K$-nilpotent objects.  It is proved in \cite[Prop.~4.9]{CI}.

\begin{thm}
\label{thm:K-nilp}
Let $K$ be a set of fibrant spaces.  In the model structure of
Theorem \ref{thm:local}, the fibrant objects are precisely the
pro-spaces that are both strictly fibrant and essentially
levelwise $K$-nilpotent.
\end{thm}

Recall that a fibration is \dfn{principal} if it is the base change
of a fibration with contractible total space.
We do not require that the base be connected.  Therefore, some of the fibers
of a principal fibration may be empty; however, the non-empty fibers are
all weakly homotopic.

\begin{lem}
\label{lem:prin}
The fibration $p:K(A,n)^{\Delta^k} \map K(A,n)^{\bd^k}$ is principal
for all $n \geq 0$, $k \geq 0$, and every abelian group $A$.  Its non-empty
fiber is weakly equivalent to $K(A, n-k)$.
\end{lem}

\begin{proof}
Throughout this proof, we use models for $K(A,n)$ that are simplicial
abelian groups \cite{Ma}.  
In particular, this means that $K(A,n)$ is fibrant.  Also,
$K(A,n)$ is based at $0$.

Only one path-component of $K(A,n)^{\bd^k}$ is in the image
of $p$; it consists of the maps $\bd^k \map K(A,n)$ that
are null-homotopic.  If $n \neq k-1$, then $K(A,n)^{\bd^k}$ only
has one component, but this is irrelevant.
Thus, we only have to compute the fiber of $p$
over the zero map $0: \bd^k \map K(A,n)$ (which is a point of
$K(A,n)^{\bd^k}$), and this fiber is
$\Omega^k K(A,n)$, as desired.

We still have to check that $p$ is principal.  Let $v$ be the $0$th
vertex of $\Delta^k$.  Then we have a short exact sequence
\[
0 \map X_{n,k} \map K(A,n)^{\Delta^k} \map K(A,n)^{\{v\}}=K(A,n) \map 0
\]
of simplicial abelian groups, where $X_{n,k}$ 
is the subspace of $K(A,n)^{\Delta^k}$
consisting of all maps that take $v$ to $0$.
The projection $\Delta^k \map \{v\}$ gives
a splitting, so $K(A,n)^{\Delta^k}$ is isomorphic to 
$K(A,n) \times X_{n,k}$.
Similarly, $K(A,n)^{\bd^k}$ is isomorphic to $K(A,n) \times Y_{n,k}$, where
$Y_{n,k}$ is the subspace of $K(A,n)^{\bd^k}$ consisting of 
all maps that take $v$ to $0$.
The fibration $X_{n,k} \map Y_{n,k}$ is principal because $X_{n,k}$
is contractible (which follows from the fact 
that $\Delta^k$ is contractible).  The identity map on $K(A,n)$ is
of course principal, and 
a product of two principal fibrations is again principal.
This shows that $p$ is principal.
\end{proof}

\begin{prop}
\label{prop:nilp-tower}
Let $C$ be a class of $R$-modules, and 
let $K$ be the class of
Eilenberg-Mac Lane spaces $K(M,n)$ such that $M$ belongs to $C$.
A space $X$ is $K$-nilpotent if and only if it is fibrant and
there exists a finite tower
\[
X_n \map X_{n-1} \map \cdots \map X_1 \map X_0 = *,
\]
where $X_n$ is weakly equivalent to $X$ and each map $X_k \map X_{k-1}$ 
is a principal fibration whose non-empty fibers belong to $K$.
\end{prop}

\begin{proof}
For one direction, Lemma \ref{lem:prin} tells us that for every $M$ in $C$,
the map $p:K(M,n)^{\Delta^k} \map K(M,n)^{\bd^k}$ is a principal fibration
whose non-empty fibers are weakly equivalent to $K(M,n-k)$.
In view of Lemma \ref{lem:K-nilp}, 
this shows that every $K$-nilpotent
space $X$ has a tower of the desired form.

Now we will show that if $X$ has a tower of the desired form, then
$X$ is $K$-nilpotent.  By induction on the length of the tower,
we just need to show that if $p:E \map B$ is a principal fibration whose
non-empty fibers are weakly equivalent to $K(M,n)$ and $B$ is $K$-nilpotent, 
then $E$ is $K$-nilpotent.  
We can rewrite $p$ in the form $\phi \coprod E \map B_0 \coprod B_1$, 
where the fibers over $B_0$ are all empty and the fibers over $B_1$
are all weakly equivalent to $K(M,n)$.  By the following lemma,
we know that $B_1$ is $K$-nilpotent.  
Thus it suffices to replace $B$ with $B_1$ and assume that $p$ is surjective.

Let $p': E' \map B'$ be a fibration with $E'$ contractible such that 
$p$ is a base change of $p'$.  Since $p$ is surjective, the image of 
the map $B \map B'$ lies in the same component as the image of $p'$.
This means that we may assume that $B'$ is connected.
Since every fiber of $p'$ is weakly equivalent to $K(M,n)$, we know
that $B'$ is weakly equivalent to $K(M,n+1)$.

We now know that $E$ is the homotopy pullback of the diagram
\[
B \map K(M,n+1) \leftarrow *.
\]
Recall the spaces $X_{n+1,1}$ and $Y_{n+1,1}$ from the proof of
Lemma \ref{lem:prin}.
Since $\bd^k$ is a model for $S^{k-1}$, the space
$Y_{n,k}$ is weakly equivalent to $\Omega^{k-1} K(M,n)$, which is
a model for $K(M,n-k+1)$.  
Thus $E$ is the homotopy pullback of the diagram
\[
B \map Y_{n+1,1} \leftarrow X_{n+1,1},
\]
so $E$ is 
also the homotopy pullback of the diagram
\[
B \map K(M,n+1) \times Y_{n+1,1} \leftarrow K(M,n+1) \times X_{n+1,1}.
\]
This identifies $E$ up to homotopy as a fiber product
\[
B \times_{K(M,n+1)^{\bd^1}} K(M,n+1)^{\Delta^1},
\]
which implies that $E$ is $K$-nilpotent.
\end{proof}

\begin{lem}
\label{lem:K-nilp-component}
Let $C$ be a class of $R$-modules, and 
let $K$ be the class of
Eilenberg-Mac Lane spaces $K(M,n)$ such that $M$ belongs to $C$.
Let $X = Y \coprod Z$.  If $X$
is $K$-nilpotent, then $Y$ and $Z$ are $K$-nilpotent.
\end{lem}

We make no assumptions about whether $Y$ and $Z$ are connected.

\begin{proof}
We choose
the model for $K(M,0)$ that consists simply of the underlying set of $M$,
viewed as a discrete simplicial set.
The map $K(M,0)^{\Delta^1} \map K(M,0)^{\bd^1}$ then becomes the
diagonal $M \map M \times M$.
Take any map $X \map K(M,0)^{\bd^1}$ such that $Y$ maps to the 
diagonal and $Z$ maps off the diagonal.  Then the pullback
$X \times_{K(M,0)^{\bd^1}} K(M,0)^{\Delta^1}$ is isomorphic to $Y$.
This shows that $Y$ is $K$-nilpotent.  The same argument shows that $Z$ is
$K$-nilpotent.
\end{proof}

When $C$ is the collection of all $R$-modules, the notion of
$K$-nilpotence defined here is related but not equivalent to
the notion of $R$-nilpotence in \cite{BK}.  However, as
Proposition \ref{prop:nilp-tower} and \cite[Prop.~III.5.3]{BK} show,
every space that is nilpotent in the sense of Bousfield and Kan is also
nilpotent in our sense.
The converse is not true.  For example, a $K$-nilpotent space 
(such as $K(R,0)$) need not be connected.  
One way to see the difference is to note that the two notions of
nilpotence are based upon different notions of principality.
The question is whether the base spaces are required to be connected
(or, equivalently, whether empty fibers are allowed).

It would be nice to have algebraic conditions on the homotopy groups
of a disconnected space that guarantee that the space is $K$-nilpotent.  
However, such a condition has eluded us so far.  See Section \ref{sctn:?}
for an elaboration of this problem.

\section{Pro-Homological Algebra}
\label{sctn:hom-alga}

In this section, we let $\A$ be any abelian category.  Recall that
the category $\pro \A$ is again abelian \cite[App.~4.5]{AM}.  
The monomorphisms
are the essentially levelwise monomorphisms, and the epimorphisms
are the essentially levelwise epimorphisms.  We also assume that
$\A$ has enough injectives.  This implies that $\pro \A$
also has enough injectives \cite{Z}.

\begin{lem}
\label{lem:constant-inj}
Let $\A$ be an abelian category, and let $I$ be an injective object of
$\A$.  Then $cI$ is an injective object of $\pro \A$.
\end{lem}

\begin{proof}
Let $i:A \map B$ be a monomorphism in $\pro \A$; we may assume that
$i$ is a level monomorphism.
Let $f: A \map cI$ be any map
in $\pro \A$.  This map is represented by a map
$f_{s}: A_{s} \map I$ in $\A$ for some $s$.  Now $f_s$ extends
over $i_s$ since $I$ is injective and $i_s$ is a monomorphism.
This extension represents a map $g: B \map I$, and $g$ extends $f$ over $i$.
\end{proof}

\begin{lem}
\label{lem:pro-ext}
Let $\A$ be an abelian category, and let $A$ be any object
of $\pro \A$.  
The groups $\Ext_{\pro \A}^{n} (A, cB)$ and 
$\colim_{t} \Ext_{\A}^{n}(A_s, B)$ are isomorphic 
for every object $B$ of $\A$ and every $n \geq 0$.
\end{lem}

\begin{proof}
Let 
\[
B \map I^0 \map I^1 \map \cdots
\]
be an injective resolution of $B$ in $\A$.  Then
\[
cB \map cI^0 \map cI^1 \map \cdots
\]
is an injective resolution of $cB$ in $\pro \A$ by 
Lemma~\ref{lem:constant-inj}
and the fact that $c$ preserves exactness.  Therefore,
$\Ext_{\pros}^{*}(A, cB)$ is the homology of the complex
\[
\Hom_{\pros} (A, cB) \map 
\Hom_{\pros} (A, cI^0) \map 
\Hom_{\pros} (A, cI^1) \map \cdots,  
\]
which is equal to the complex
\[
\colim_s \Hom_{\A} (A_s, B) \map
\colim_s \Hom_{\A} (A_s, I^0) \map
\colim_s \Hom_{\A} (A_s, I^1) \map \cdots.
\]
Since filtered colimits are exact, the homology of the last complex is
equal to $\colim_s \Ext_{\A}^{*} (A_s, B)$.
\end{proof}

Now we construct a universal coefficients spectral sequence
for pro-spaces.  Let $X$ be an arbitrary pro-space, and let $M$
be an $R$-module.  Choose an injective
resolution 
\[
0 \map M \map I^0 \map I^1 \map \cdots.
\]
Consider the bicomplex $K^{*,*}$ given by the formula 
\[
K^{p,q} = C^{q} (X; I^{p}) = \colim_s C^{q} (X_s; I^{p})
= \colim_s \Hom(C_q (X_s;R), I^p).
\]
This is a first-quadrant bicomplex with cohomological grading.
Now $C_q(X_s;R)$ is 
a free $R$-module, so the complex 
$\colim_s \Hom(C_q(X_s;R), I^*)$ is exact.  Taking cohomology of
$K^{*,*}$ with respect to the $p$-differential gives that
\[
E_1^{0,q} =
\colim_s \Hom(C_q(X_s;R), M) = C^q(X; M)
\]
and $E_1^{p,q} = 0$ if $p > 0$.  Therefore,
$E_2^{0,q} = H^q(X;M)$ and $E_2^{p,q} = 0$ if $p > 0$.
Thus, the spectral sequence collapses, and the
cohomology of the total complex of $K^{*,*}$ is $H^*(X;M)$.

Now we compute in the other order.  
Taking cohomology of $K^{*,*}$ with respect to the $q$-differential gives
\[
E_{1}^{p,q} = H^{q} (X; I^{p}) = \colim_s H^{q} (X_s; I^{p}).
\]
Because $I^p$ is injective, this equals
$\colim_s \Hom( H_q(X_s;R), I^p)$.  Therefore,
taking cohomology with respect to the $p$-differential gives
\[
E_{2}^{p,q} = 
\colim_s \Ext\nolimits_R^p(H_q(X_s;R), M) =
\Ext\nolimits_{\pros}^p (H_q(X;R), M).
\]
The second equality above relies on Lemma \ref{lem:pro-ext}.

Hence, we have a convergent first-quadrant cohomological spectral 
sequence
\[
E_{2}^{p,q} = \Ext\nolimits_{\pros}^{p} ( H_{q}(X; R), cM) \Longrightarrow
  H^{p+q}(X;M).
\]
This spectral sequence is called the \dfn{pro-universal coefficients
spectral sequence}.

\section{Cohomological and Homological Weak Equivalences}
\label{sctn:wk-eq}

In this section we collect some results about the cohomology and pro-homology
of pro-spaces, emphasizing the differences and similarities
with the situation of ordinary spaces.

The first lemma looks quite strange at first glance, but it becomes
plausible when one remembers that filtered limits are exact in pro-categories
\cite{I1} \cite{AHJM}.

\begin{lem}
\label{lem:lim-cohlgy}
If $a \mapsto X^a$ is a cofiltered diagram of pro-spaces, then
the cohomology group $H^n(\lim_a^{\pros} X^a; M)$ is isomorphic to
$\colim_a H^n(X^a; M)$, where $\lim^{\pros}$ is the limit internal to
the category of pro-spaces.
\end{lem}

\begin{proof}
This follows by direct computation using the construction of cofiltered limits
in pro-categories given in Section \ref{subsctn:pro-ctgy}.
\end{proof}

\begin{defn}
\label{defn:R-cohlgy}
Let $R$ be any ring.
An \mdfn{$R$-cohomology weak equivalence} 
is a map of pro-spaces $X \map Y$
inducing isomorphisms
$H^{n} (Y; R) \map H^{n} (X; R)$
for every $n \geq 0$.
\end{defn}

Suppose that $f$ is a map of ordinary spaces inducing an $R$-cohomology
isomorphism.  Then
$f$ induces an isomorphism in
cohomology with coefficients in an arbitrary product of copies of $R$.
But every free $R$-module is a retract of some product
of copies of $R$,
so $f$ induces an isomorphism in cohomology with coefficients
in any free $R$-module.  Now every $R$-module $M$ belongs to
a short exact sequence
\[
0 \map F_1 \map F_2 \map M \map 0
\]
in which $F_1$ and $F_2$ are free, so the long exact sequence of 
cohomology groups and the five lemma imply that $f$ induces an
isomorphism in cohomology with coefficients in any $R$-module.

In contrast to the above paragraph, if $f$ is a map of pro-spaces that
is an $R$-cohomology weak equivalence, then $f$ does not 
necessarily induce a cohomology isomorphism with coefficients in 
all $R$-modules.
The argument from the previous paragraph 
breaks down because cohomology of pro-spaces
does not commute with arbitrary products
of coefficients.  Actually, cohomology only commutes with finite
products.  One explanation for this difference is that
$K( \prod_{t} M_{t}, n)$ is weakly equivalent to
$\prod_{t} K(M_{t}, n)$ for ordinary spaces, but 
$cK(\prod_{t} M_{t}, n )$ is not weakly equivalent to
$\prod_{t} cK( M_{t}, n)$ 
as pro-spaces.

\begin{ex}
\label{ex}
Let $V$ be an infinite-dimensional $\Z/p$-vector space, and let $n \geq 1$.
We give an example of a pro-space $X$ for which
$H^n(X; \Z/p) = 0$ but $H^n(X; V)$ is non-zero.
Consider the pro-vector space $W$ consisting of all subspaces of $V$
with finite codimension; the structure maps are the inclusions
of subspaces.  Let $X$ be the pro-space $K(W, n)$ obtained by applying
levelwise the functor $K(-,n)$ to $W$.

Now $H^n(X; \Z/p)$ equals $\colim_s \Hom_{\Z/p}(W_s, \Z/p)$.
This colimit is zero because the kernel of any homomorphism
$W_s \map \Z/p$ is equal to $W_t$ for some $t$.
On the other hand, $H^n(X; V)$ equals $\colim_s \Hom_{\Z/p}(W_s, V)$,
which is non-zero because it contains the element represented
by any of the inclusions $W_s \map V$.
\end{ex}

Because of this phenomenon, we introduce the following definition,
which is distinct from Definition \ref{defn:R-cohlgy}.

\begin{defn}
\label{defn:R-hlgy}
Let $R$ be any ring.
An \mdfn{$R$-homology weak equivalence} 
is a map of pro-spaces $X \map Y$ inducing isomorphisms 
$H_n (X; R) \map H_n (Y; R)$
of pro-groups for all $n \geq 0$.
\end{defn}

This definition can be reformulated in terms of cohomology.  
The next result implies that every $R$-homology weak equivalence
is an $R$-cohomology weak equivalence.
As shown in Example \ref{ex}, the converse is not true.

\begin{prop}
\label{prop:hlgy-eq}
A map $f:X \map Y$ is an $R$-homology weak equivalence if and only if
it induces an isomorphism $f^*: H^{n} (Y; M) \map H^{n}(X; M)$
for all $n \geq 0$ and all $R$-modules $M$.
\end{prop}

The following proof is inspired by \cite[Prop.~III.6.7]{BK}.

\begin{proof}
First suppose that $f$ is an $R$-homology weak equivalence.
For any $R$-module $M$, $f$ induces an 
isomorphism of $E_2$-terms of the pro-universal
coefficients spectral sequence (see Section \ref{sctn:hom-alga}).  
Therefore, the map on abutments is also an isomorphism.

For the other implication,
suppose that $f$ is a cohomology isomorphism with coefficients
in any $R$-module.  
If $I$ is an injective $R$-module, then
\[
H^{n} (X; I) =
\colim_{s} \Hom ( H_{n} (X_{s}; R), I) =
\Hom_{\pros}(H_{n} (X; R), cI).
\]
Similarly, $H^{n} (Y; I) = \Hom_{\pros}(H_{n} (Y; R), cI)$. 
It follows that
\[
\Hom_{\pros}(H_{n} (Y; R), cI) \map
\Hom_{\pros}(H_{n} (X; R), cI) 
\]
is an isomorphism for every injective $R$-module $I$.

Now let $M$ be an arbitrary pro-$R$-module.  Choose a monomorphism
$M \map I$ such that $I$ is an injective pro-$R$-module.
Then there is a diagram
\[
\xymatrix{
\Hom_{\pros} (H_{n} (Y; R), M) \ar[r] \ar[d] &
  \Hom_{\pros} (H_{n} (Y; R), I) \ar[d] \\
\Hom_{\pros} (H_{n} (X; R), M) \ar[r] &
  \Hom_{\pros} (H_{n} (X; R), I)            }
\]
in which the horizontal maps are monomorphisms.  
Since the right vertical map is
an isomorphism, the left vertical map must be a monomorphism.
We conclude that $H_n(X;R) \map H_n(Y;R)$ is an epimorphism of
pro-abelian groups.

Now let $K$ be the pro-group that is the kernel of the map
$H_n(X;R) \map H_n(Y;R)$.  Since $\Hom_{\pros}(-,I)$ is exact for
all injective pro-groups $I$,
the first paragraph tells us that 
$\Hom_{\pros}(K,I)$ is zero for all injectives $I$.  
But the category of
pro-abelian groups has enough injectives, so this can only happen
if $K$ equals zero.  This means that the map
$H_n(X;R) \map H_n(Y;R)$ is an isomorphism.
\end{proof}

From now on, we will freely switch between the cohomological
and homological descriptions of $R$-homology weak equivalences,
as given in Definition \ref{defn:R-hlgy} and 
Proposition \ref{prop:hlgy-eq}.

We will need to rephrase cohomology isomorphisms in terms of
weak equivalences of mapping spaces.

\begin{prop}
\label{prop:cohlgy-wk-eq}
Let $M$ be any $R$-module, 
and let $K(M,n)$ be a fibrant Eilenberg-Mac Lane space.
If $f: X \map Y$ is any map of pro-spaces,
then $H^n(f;M)$ is an isomorphism for all $n \geq 0$ if and only if the map
$\Map_{\pros}(Y, cK(M,n)) \map \Map_{\pros}(X, cK(M,n))$ is a weak equivalence
for all $n \geq 0$.
\end{prop}

\begin{proof}
First suppose that the maps between mapping spaces are weak equivalences.
We get isomorphisms after taking $\pi_0$, which gives us the desired
cohomology isomorphisms.

Now suppose that the maps $H^n(f;M)$ are isomorphisms for all $n\geq 0$.  
Since all pro-spaces
are cofibrant, the homotopy types of the mapping spaces do not change if
we alter $X$ or $Y$ up to strict weak equivalence.  This means that we may
assume that $f$ is a levelwise cofibration; let $Z$ be its cofiber,
which is computed levelwise.  Recall that $Z$ is pointed canonically.
Using the long exact sequence in cohomology, note that the
reduced cohomology group $\tilde{H}^n(Z;M)$
is zero for all $n \geq 0$.

Now we want to show that the fibration 
\[
\Map_{\pros}(Y, cK(M,n)) \map \Map_{\pros}(X, cK(M,n))
\]
is actually an acyclic fibration of simplicial sets.  
We can do this by showing that it has
the right lifting property with respect to all generating cofibrations
$\bd^k \map \Delta^k$.  After the usual adjointness arguments, we need
to show that every diagram
\[
\xymatrix{
\bd^k \otimes Y \coprod_{\bd^k \otimes X} \Delta^k \otimes X \ar[r] \ar[d] &
cK(M,n) \\
\Delta^k \otimes Y       }
\]
has a lift.  Since $K(M,n)$ is fibrant, formal model category arguments
show that we only need obtain a lift up to homotopy; then we can adjust
this lift to get an actual lift.
In other words, we need to show that the map
\[
\xymatrix@1{
H^n(\Delta^k \otimes Y; M) \ar[r] &
  H^n(\bd^k \otimes Y \coprod_{\bd^k \otimes X} \Delta^k \otimes X; M)  }
\]
is surjective.
The cofiber of the vertical map
above is $S^k \smash Z$, where the smash product is constructed levelwise.
Using the long exact sequence in cohomology,
it suffices to show that the reduced cohomology group
$\tilde{H}^{n+1}(S^k \smash Z; M)$ is zero.  But this group
is equal to $\tilde{H}^{n+1-k}(Z;M)$, which we already computed to be zero.
\end{proof}

From now on, we will frequently express cohomology isomorphisms in terms
of weak equivalences of mapping spaces as in 
Proposition \ref{prop:cohlgy-wk-eq}.
For example, we have the following corollary.

\begin{cor}
\label{cor:cohlgy-wk-eq}
\mbox{}
\begin{enumerate}
\item A map $f$ is an $R$-cohomology weak equivalence if and only if
the map 
\[
\Map_{\pros}(f, cK(R,n))
\]
is a weak equivalence for every $n \geq 0$, where $K(R,n)$ is a fibrant
Eilenberg-Mac Lane space.
\item A map $f$ is an $R$-homology weak equivalence if and only if
the map 
\[
\Map_{\pros}(f, cK(M,n))
\]
is a weak equivalence for every $n \geq 0$ and 
every $R$-module $M$,
where $K(M,n)$ is a fibrant Eilenberg-Mac Lane space.
\end{enumerate}
\end{cor}

\begin{proof}
The first claim follows immediately from 
Definition \ref{defn:R-cohlgy} and Proposition \ref{prop:cohlgy-wk-eq}.
The second claim follows immediately from
Propositions \ref{prop:hlgy-eq} and \ref{prop:cohlgy-wk-eq}.
\end{proof}

\section{The $R$-Cohomological and $R$-Homological Model Structures}
\label{sctn:cohlgy-mc}

We next establish model structures whose weak equivalences are 
the $R$-cohom\-ology weak equivalences and the $R$-homology weak equivalences.

\begin{defn}
A \dfn{cofibration} of pro-spaces is an essentially levelwise cofibration.
\end{defn}

This is the same notion of cofibration as in the 
strict model structure \cite{I2} or the $\pi_*$-model structure \cite{I3}.
In fact, the notion of cofibration is the same for every model
structure on pro-spaces in the entire paper.

\begin{defn}
An \mdfn{$R$-cohomology fibration} 
is a map having the right lifting property
with respect to all maps that are both cofibrations
and $R$-cohomology weak equivalences.
\end{defn}

\begin{thm}
\label{thm:cohlgy-mc}
The cofibrations, $R$-cohomology weak equivalences, and 
$R$-cohom\-ology fibrations give a left proper simplicial model structure
on the category of pro-spaces.
\end{thm}

\begin{proof}
The proof is an application of Theorem \ref{thm:local}
with $K$ equal to the 
set of Eilenberg-Mac Lane objects of the form $K(R, n)$.
By Corollary \ref{cor:cohlgy-wk-eq}(1), the
$K$-colocal weak equivalences are the same as the 
$R$-cohomology weak equivalences.
\end{proof}

If $\lambda$ is any infinite cardinal, then an $R$-module $M$ is
\mdfn{$\lambda$-generated} if $M$ has a generating set of cardinality at most
$\lambda$.  Note that there is a set of isomorphism types of 
$\lambda$-generated $R$-modules.  
Also note that a $\lambda$-generated $R$-module might have more than
$\lambda$ elements because $R$ might be large.

A pro-map is a \mdfn{$\lambda$-generated $R$-cohomology weak equivalence}
if it induces an isomorphism on cohomology with coefficients in
all $\lambda$-generated $R$-modules.  Similarly,
a pro-map is a \mdfn{$\lambda$-generated $R$-cohomology fibration}
if it has the right lifting property with respect to all pro-maps that
are both cofibrations and $\lambda$-generated $R$-cohomology weak 
equivalences.

\begin{thm}
\label{thm:lambda-cohlgy-mc}
For any infinite cardinal $\lambda$,
the cofibrations, $\lambda$-generated $R$-co\-ho\-mology 
weak equivalences, and 
$\lambda$-generated $R$-cohomology fibrations 
give a left proper simplicial model structure
on the category of pro-spaces.
\end{thm}

\begin{proof}
The proof is an application of Theorem \ref{thm:local}
with $K$ equal to the 
set of Eilenberg-Mac Lane objects of the form $K(M, n)$ with
$M$ a $\lambda$-generated $R$-module.
\end{proof}

Actually, the $\lambda$-generated model structures are not so interesting.
What we really want is a model structure in which the weak equivalences
are detected by $R$-homology, or equivalently, 
according to Proposition \ref{prop:hlgy-eq},
by cohomology with coefficients in all $R$-modules.  
This is the main purpose of the rest of this section.

\begin{defn}
An \mdfn{$R$-homology fibration} 
is any map that is a $\lambda$-generated $R$-cohomology fibration
for some $\lambda$.
\end{defn}

In other words, the class of $R$-homology fibrations is the union of the
classes of $\lambda$-generated $R$-cohomology fibrations as $\lambda$
ranges over all cardinals.
By Proposition \ref{prop:hlgy-eq}, 
the class of $R$-homology weak equivalences is the 
intersection of the classes of $\lambda$-generated $R$-cohomology weak
equivalences as $\lambda$ ranges over all cardinals.

\begin{lem}
\label{lem:lambda-fib}
For any $\lambda$,
the acyclic $\lambda$-generated $R$-cohomology fibrations  
are the same as the acyclic $R$-homology fibrations.
\end{lem}

\begin{proof}
For every $\lambda$, the acyclic $\lambda$-generated $R$-cohomology 
fibrations are detected by the same class of cofibrations, so these
classes of acyclic $\lambda$-generated $R$-cohomology fibrations
are all equal.
Therefore,
if $p$ is an acyclic $\lambda$-generated
$R$-cohomology fibration, then it is a
$\mu$-generated $R$-cohomology weak equivalence for all $\mu$.
This implies that $p$ is an $R$-homology weak equivalence
by Proposition \ref{prop:hlgy-eq}.

For the other direction, suppose that $p$ is an acyclic $R$-homology
fibration.  
Then it is a $\mu$-generated $R$-cohomology weak equivalence
for every $\mu$ and a $\lambda$-generated $R$-cohomology fibration for
some $\lambda$.  Thus, $p$ is an acyclic $\lambda$-generated $R$-cohomology
fibration for some $\lambda$.  As in the first paragraph, 
this implies that $p$ is an acyclic $\lambda$-generated $R$-cohomology
fibration for all $\lambda$.
\end{proof}

\begin{thm}
\label{thm:hlgy-mc}
The cofibrations, $R$-homology weak equivalences, and 
$R$-homology fibrations are a left proper simplicial model structure
on the category of pro-spaces.
\end{thm}

\begin{proof}
Most of the proof follows formally from the existence of the 
$\lambda$-generated $R$-cohomology model structures.  Given any
collection of classes that satisfy the two-out-of-three axiom,
their intersection also satisfies the two-out-of-three axiom.
Similarly, arbitrary unions and intersections of classes preserve
the retract axiom.

Since Lemma \ref{lem:lambda-fib} tells us what
the acyclic $R$-homology fibrations are,
one of the lifting axioms and one of the factoring axioms follows immediately
from Theorem \ref{thm:lambda-cohlgy-mc}.  
For the other lifting axiom, suppose that $i$
is an acyclic $R$-homology cofibration and $p$ is an $R$-homology fibration.
Then there exists some $\lambda$ such that $p$ is a $\lambda$-generated
$R$-cohomology fibration, and $i$ is also an acyclic $\lambda$-generated
$R$-cohomology cofibration.  Thus $i$ has the left lifting property
with respect to $p$ because of Theorem \ref{thm:lambda-cohlgy-mc}.

Factorizations into acyclic $R$-homology cofibrations followed by
$R$-homology fibrations are significantly more difficult.  We construct
these below in Proposition \ref{prop:factor-hlgy-mc}.

The simplicial structure also follows formally from the
$\lambda$-generated model structures of Theorem
\ref{thm:lambda-cohlgy-mc}.  Namely, let $i: A \map B$
be a cofibration and let $p: X \map Y$ be an $R$-homology fibration.
Then $p$ is a $\lambda$-generated $R$-homology fibration for some $\lambda$,
so the map
\[
\Map(i,p): \Map(B, X) \map \Map(A, X) \times_{\Map(A, Y)} \Map(B, Y)
\]
is a fibration of simplicial sets because the $\lambda$-generated
model structure is simplicial.  If $i$ or $p$ is acyclic, then it
is a $\lambda$-generated $R$-cohomology weak equivalence, so the 
above map $\Map(i,p)$ is also a weak equivalence because
of the $\lambda$-generated model structure of 
Theorem \ref{thm:lambda-cohlgy-mc}.

Left properness follows immediately from the fact that every object
is cofibrant. 
\end{proof}

Functorially in $M$, 
choose a fibrant Eilenberg-Mac Lane space $K(M,n)$ for each finitely
generated $R$-module $M$.  Let $p^k(M,n)$ be the fibration
$K(M,n)^{\Delta^k} \map K(M,n)^{\bd^k}$.

Now for any $R$-module $M$, define $K(M,n)$ to be
$\colim K(N,n)$,
where the colimit ranges over all finitely generated submodules $N$
of $M$.  This construction is a specific fibrant model for the 
Eilenberg-Mac Lane space of type $(M,n)$.
We also get maps $p^k(M,n)$, which are again fibrations
because the colimit is filtered.  Note also that $K(M,n)^{\Delta^k}$
is equal to $\colim_{N \leq M} (K(N,n)^{\Delta^k})$ because
$\Delta^k$ is a finite simplicial set (and similarly for
$K(M,n)^{\bd^k}$).

The constant map
$cp^k(M,n): cK(M,n)^{\Delta^k} \map cK(M,n)^{\bd^k}$ 
is a $\lambda$-generated $R$-cohomology fibration 
if $M$ is $\lambda$-generated.  This (or rather its dual)
is proved in \cite[Lem.~2.3(b)]{CI}; it can also
be deduced from Proposition \ref{prop:cohlgy-wk-eq} using adjointness.
Therefore, $cp^k(M,n)$ is an $R$-homology fibration for every $M$.

A pro-space $X$ is \mdfn{$\lambda$-bounded}
if each space $X_s$ has at most $\lambda$ elements.

\begin{lem}
\label{lem:lambda-lift}
Let $f:X \map Y$ be any map between $\lambda$-bounded pro-spaces.
Then $f$ has the left lifting property with respect to all maps
of the form $cp^k(M,n)$ with $M$ any $R$-module, if and only if 
$f$ has the left lifting property with respect to all maps
of the form $cp^k(M,n)$ with $M$ any $\lambda$-generated $R$-module.
\end{lem}

\begin{proof}
One direction is tautological.  For the other direction, 
let $f$ have the left lifting property with respect to all maps
of the form $cp^k(M,n)$ with $M$ any $\lambda$-generated $R$-module.
Suppose there is a square
\[
\xymatrix{
X \ar[r]\ar[d] & cK(M,n)^{\Delta^k} \ar[d] \\
Y \ar[r] & cK(M,n)^{\bd^k}                 }
\]
with $M$ any $R$-module; we want to produce a lift.
This square can be represented as a square
\[
\xymatrix{
X_s \ar[r]^-g \ar[d] & K(M,n)^{\Delta^k} \ar[d] \\
Y_t \ar[r]_-h        & K(M,n)^{\bd^k}        }
\]
of ordinary spaces.
For each $x$ in $X_s$, $f(x)$ belongs to $K(P_x, n)^{\Delta^k}$ for some
finitely generated submodule $P_x$ of $M$.  Similarly,
for each $y$ in $Y_t$, $g(y)$ lies in $K(P_y, n)^{\bd^k}$ for some
finitely generated submodule $P_y$ of $M$.  Since there
are at most $\lambda$ choices for $x$ and $y$, we can choose a
$\lambda$-generated submodule $N$ containing each $P_x$ and each $P_y$.
Now we have a diagram
\[
\xymatrix{
X \ar[r]\ar[d] & cK(N,n)^{\Delta^k} \ar[r]\ar[d] & 
   cK(M,n)^{\Delta^k} \ar[d] \\
Y \ar[r] & cK(N,n)^{\bd^k} \ar[r] & cK(M,n)^{\bd^k}                 }
\]
for some $\lambda$-generated submodule $N$ of $M$.  A lift exists in
the left square by assumption, and this gives us the desired lift.
\end{proof}

The importance of Lemma \ref{lem:lambda-lift} 
is that the second condition
involves only a set of maps of the form $cp^k(M,n)$, while the first 
condition does not.

\begin{lem}
\label{lem:lift-hlgy-mc}
A cofibration $i: A \map B$ induces an isomorphism in cohomology with
coefficients in $M$ if and only if it has the right lifting property
with respect to the maps $cp^k(M,n)$.
\end{lem}

\begin{proof}
By Proposition \ref{prop:cohlgy-wk-eq}, instead of considering
whether $H^n(i;M)$ is an isomorphism for all $n$, we shall consider whether
\[
\Map_{\pros}(i,cK(M,n)): \Map_{\pros}(B,cK(M,n)) \map \Map_{\pros}(A,cK(M,n))
\]
is a weak equivalence for all $n$.
Since $i$ is a cofibration, this map is a fibration of simplicial sets.
It is an acyclic fibration if and only if it has the left lifting property
with respect to the maps $\bd^k \map \Delta^k$.  By adjointness,
this lifting property is equivalent to the desired lifting property.
\end{proof}

\begin{prop}
\label{prop:factor-hlgy-mc}
Any map $f:X \map Y$ of pro-spaces factors into 
$i:X \map Z$ followed by $q:Z \map Y$, where
$i$ is an acyclic $R$-homology cofibration and $q$ is an $R$-homology
fibration.
\end{prop}

\begin{proof}
Choose a cardinal $\lambda_0$ such that $X$ and $Y$ are both
$\lambda_0$-bounded, and let $f_0 = f$.  Factor the map $f_0$ into
an acyclic $\lambda_0$-generated $R$-cohomology cofibration 
$f_1: X \map Z_1$ followed
by a $\lambda_0$-generated $R$-cohomology fibration $Z_1 \map Y$.  
Here we are using
Theorem \ref{thm:lambda-cohlgy-mc}.

Proceeding inductively,
choose a cardinal $\lambda_n$ such that $X$ and $Z_n$ are both
$\lambda_n$-bounded.  Factor the map $f_{n}$ into
an acyclic $\lambda_n$-generated cofibration $f_{n+1}: X \map Z_{n+1}$ followed
by a $\lambda_n$-generated fibration $Z_{n+1} \map Z_n$. 
This process yields a diagram
\[
X \map \cdots \map Z_2 \map Z_1 \map Y,
\]
and we let $Z = \lim^{\pros}_n Z_n$, where the limit is computed within
the category of pro-spaces.

Choose a cardinal $\lambda_\infty$ such that $\lambda_\infty \geq \lambda_n$
for each $n$.  Note that $X$, $Y$, and each $Z_n$ are 
$\lambda_\infty$-bounded.  Using the explicit construction of cofiltered
limits in pro-categories given in Section \ref{subsctn:pro-ctgy}, 
it follows that $Z$
is also $\lambda_\infty$-bounded.

If $\lambda \leq \mu$, then the $\mu$-generated $R$-cohomology weak
equivalences are contained in the $\lambda$-generated $R$-cohomology 
weak equivalences.  Therefore, the acyclic $\mu$-generated
$R$-cohomology cofibrations are contained in the 
$\lambda$-generated $R$-co\-hom\-ology cofibrations.  From the nature
of lifting properties, it follows that the $\lambda$-generated
$R$-cohomology fibrations are contained in the $\mu$-generated
$R$-cohomology fibrations.

Since $\lambda_n \leq \lambda_\infty$, 
each map $Z_n \map Z_{n-1}$ (and also $Z_1 \map Y$)
is a $\lambda_\infty$-generated $R$-cohomology fibration.  Thus, the
map $q:Z \map Y$ is a composition of a countable tower of
$\lambda_\infty$-generated $R$-cohomology fibrations.
Formal arguments with lifting properties imply that $q$
is also
a $\lambda_\infty$-generated $R$-cohomology fibration.  This means
that $q$ is an $R$-homology fibration.

Each map $X \map Z_n$ is a cofibration.  Since cofibrations are
closed under cofiltered limits \cite[Cor.~5.3]{I1}, we conclude that 
$i:X \map Z$ is a cofibration.  By Lemma \ref{lem:lift-hlgy-mc},
it remains only to show that $i$ has the left lifting property 
with respect to the maps $cp^k(M,n)$ for all $R$-modules $M$.

Suppose given a square
\[
\xymatrix{
X \ar[r]\ar[d] & cK(M,n)^{\Delta^k} \ar[d] \\
Z \ar[r] & cK(M,n)^{\bd^k}                }
\]
with $M$ a $\lambda_\infty$-generated $R$-module.
Using the explicit construction of filtered limits of pro-objects 
given in Section \ref{subsctn:pro-ctgy}, this diagram factors as
\[
\xymatrix{
X \ar[rr]\ar[d] & & cK(M,n)^{\Delta^k} \ar[d] \\
Z \ar[r] & Z_j \ar[r] & cK(M,n)^{\bd^k}                }
\]
for some $j$.  Because $X$ and $Z_j$ are both $\lambda_j$-bounded,
this diagram further factors into
\[
\xymatrix{
X \ar[r]\ar[d] & cK(N,n)^{\Delta^k} \ar[r]\ar[d] & cK(M,n)^{\Delta^k} \ar[d] \\
Z_j \ar[r] & cK(N,n)^{\bd^k} \ar[r] & cK(M,n)^{\bd^k}                }
\]
as in the proof of Lemma \ref{lem:lambda-lift}, where 
$N$ is a $\lambda_j$-generated submodule of $M$.

Lemma \ref{lem:lift-hlgy-mc} tells us that 
there is a lift in the diagram
\[
\xymatrix{
X \ar[rr]\ar[d] & & cK(N,n)^{\Delta^k} \ar[r]\ar[d] & 
       cK(M,n)^{\Delta^k} \ar[d] \\
Z_{j+1} \ar[r] \ar@{.>}[urr]\ar[r] & Z_j \ar[r] & 
   cK(N,n)^{\bd^k} \ar[r] &  cK(M,n)^{\bd^k}     }
\]
because $X \map Z_{j+1}$ is an acyclic $\lambda_{j}$-generated cofibration
and $N$ is a $\lambda_j$-generated $R$-module.
This gives us the desired lift.
\end{proof}

\section{$R$-Completions}
\label{sctn:R-comp}

The model structures of Theorems \ref{thm:cohlgy-mc} and \ref{thm:hlgy-mc}
allow us to define the $R$-completion of any pro-space.

\begin{defn}
\label{defn:completion}
Let $X$ be a pro-space.
The \mdfn{cohomological $R$-completion $\ccomp{X}$ of $X$}
is a fibrant replacement for $X$ in the cohomological model structure
of Theorem \ref{thm:cohlgy-mc}.
The \mdfn{homological $R$-completion $\hcomp{X}$ of $X$}
is a fibrant replacement for $X$ in the homological model structure
of Theorem \ref{thm:hlgy-mc}.
\end{defn}

Philosophically, the cohomological $R$-completion of a pro-space $X$
should preserve the cohomology of $X$ with coefficients in $R$
but forget all other information.
Similarly,
the homological $R$-completion of $X$
should preserve the cohomology of $X$ with coefficients in any $R$-module
but forget all other information.  
Thus, homological $R$-completion contains more information than
the cohomological $R$-completion.
See Example \ref{ex} for a pro-space $X$ such that $\hcomp{X}$ and $\ccomp{X}$
are distinct.

These ideas are made precise in the following theorem.

\begin{thm}
\label{thm:completion}
A map $f: X \map Y$ 
of pro-spaces (or of ordinary spaces) induces an isomorphism
in cohomology with coefficients in $R$ if and only if
$\ccomp{X}$ and $\ccomp{Y}$ are simplicially homotopy equivalent pro-spaces.
Similarly, the map $f$ induces an isomorphism
in cohomology with coefficients in all $R$-modules if and only if
$\hcomp{X}$ and $\hcomp{Y}$ are simplicially homotopy equivalent pro-spaces.
\end{thm}

Of course, if two pro-spaces are simplicially homotopy equivalent, then
they are strictly weakly equivalent or
equivalent with respect to any reasonable notion of homotopy theory
for pro-spaces.

\begin{proof}
The maps $X \map \ccomp{X}$ 
and $Y \map \ccomp{Y}$ are
$R$-cohomological weak equivalences.
Thus $H^*(f;R)$
is an isomorphism if and only if
$\ccomp{X}$ and $\ccomp{Y}$ are $R$-cohomologically weakly equivalent.
Since every pro-space is cofibrant, 
$\ccomp{X}$ and $\ccomp{Y}$ are both cofibrant and fibrant with respect to the
$R$-cohomological model structure of Theorem \ref{thm:cohlgy-mc}.  
Thus, $\ccomp{X}$ and $\ccomp{Y}$ are $R$-cohomologically weakly equivalent
if and only if they are simplicially homotopy equivalent.

The argument for homological $R$-completions is identical except that 
it uses the $R$-homological model structure of Theorem \ref{thm:hlgy-mc}.
\end{proof}

For any pro-space $X$ such that each $X_s$ is pointed and connected,
we will now show how to construct $\hcomp{X}$ in terms of the
Bousfield-Kan $R$-towers \cite{BK} of the spaces $X_s$.  The moral is that at 
least for pointed connected spaces, 
the Bousfield-Kan $R$-tower (with a minor modification) is the same thing as 
fibrant replacement in the $R$-homological model structure.

\begin{prop}
\label{prop:BK-tower}
Let $X$ be a pro-object in the category of pointed connected spaces, and
let $I$ be the indexing category of $X$.
Construct a new pro-space $Y$ with indexing category $I \times \N$
by defining $Y_{s,n}$ to be the $n$th Postnikov section $P_n R_n X_s$
of the $n$th stage of the Bousfield-Kan $R$-tower for $X_s$.
Then the strict fibrant replacement $\hat{Y}$ of $Y$ is a 
fibrant replacement for $X$ in the $R$-homological
model structure.
\end{prop}

\begin{rem}
In the previous proposition, it is also possible to define
$Y_{s,n}$ to be just $R_n X_s$, not its $n$th Postnikov section.  However,
then $\hat{Y}$ must be a $\pi_*$-fibrant replacement, i.e., a fibrant
replacement in the model structure on pro-spaces in which weak equivalences
are detected by pro-homotopy groups \cite{I2}.
\end{rem}

\begin{proof}
Let $K$ be the collection
of all fibrant spaces of the form $K(M,n)$, where $n \geq 0$ and $M$ is
any $R$-module.
From \cite[Cor.~III.5.6]{BK} and \cite[Prop.~III.5.3]{BK}, we know that
the Postnikov tower of $R_n X_s$ can be refined to a sequence of 
principal fibrations whose fibers belong to $K$.  Therefore, the
Postnikov tower of $P_n R_n X_s$ consists of a finite sequence
of principal fibrations whose fibers belong to $K$.
Proposition \ref{prop:nilp-tower} tells us that each 
$P_n R_n X_s$ is $K$-nilpotent, so $Y$ is essentially levelwise 
$K$-nilpotent.
It follows from \cite[Prop.~3.7]{CI} that $\hat{Y}$ is also essentially
levelwise $K$-nilpotent, so Theorem \ref{thm:K-nilp} tells us that $\hat{Y}$
is fibrant in the $R$-homological model structure.

It remains to show that the map $X \map \hat{Y}$ is an $R$-homology
weak equivalence.  Since $Y \map \hat{Y}$ is a strict weak equivalence, it 
suffices to show that $X \map Y$ is an $R$-homology weak equivalence.
We know from \cite[Prop.~III.6.5]{BK} (or \cite{D})
that $cH_k(X_s;R) \map H_k(P_* R_* X_s;R)$ is 
a pro-isomorphism for each $s$; here $P_* R_* X_s$ is the pro-space
\[
\cdots \map P_2 R_2 X_s \map P_1 R_1 X_s \map P_0 R_0 X_s,
\]
and $cH_k(X_s;R)$ is the constant pro-group with value $H_k(X_s;R)$.
Now $H_k(X;R)$ is isomorphic to $\lim^{\pros}_s cH_k(X_s;R)$, where
the limit is computed within the category of pro-groups.  Similarly,
$H_k(Y;R)$ is isomorphic to $\lim^{\pros}_s H_k(P_* R_* X_s;R)$
(see the construction of limits in pro-categories given in 
Section \ref{subsctn:pro-ctgy}).
Thus, the map $H_k(X;R) \map H_k(Y;R)$ is a pro-isomorphism.
\end{proof}

\section{$\Z/p$-Cohomology}
\label{sctn:Morel}

In this section, fix a prime $p$, and let $R$ be the finite ring
$\Z/p$.  We will compare the $\Z/p$-cohomological model structure of
Theorem \ref{thm:cohlgy-mc} to the $\Z/p$-cohomological model
structure of \cite{Mo} and show that they yield the same homotopy
categories.  Throughout this section, whenever we
discuss the category of pro-simplicial sets, we are always thinking of it
equipped with the $\Z/p$-cohomological model structure.

We first recall some ideas from \cite{Mo}.  Let $\F$ be the 
category of finite sets.  The category $\pro \F$ of pro-finite
sets is equivalent to the category of totally disconnected compact Hausdorff
topological spaces.  

The main object of study in \cite{Mo} is the category $s \pro \F$ of
simplicial pro-finite sets (or, equivalently, simplicial totally
disconnected compact Hausdorff topological spaces).  
The weak equivalences in this category are the continuous cohomology
isomorphisms with $\Z/p$-coefficients, the cofibrations are the degreewise
monomorphisms, and the fibrations are defined by a lifting property.

The main purpose of this model structure on $s \pro \F$ is to describe
a $\Z/p$-completion functor.  Given any set $X$, let $\hat{X}$ be
the pro-finite set of all finite quotients of $X$.  
Applying this construction degreewise gives a functor from
simplicial sets to simplicial pro-finite sets.
The $\Z/p$-completion of a simplicial set $X$ is defined to be a fibrant
replacement for $\hat{X}$.  

In order to compare the category $s \pro \F$ to the category of
pro-simplicial sets, we need the intermediate category $\pro s\F$ of
pro-simplicial finite sets.  
Despite claims 
in \cite{Mo}, \cite{R}, and elsewhere, this category
is not equivalent to $s \pro \F$.  See \cite[Ex.~3.7]{I1} for a counterexample.
Beware that a simplicial finite set is not the same as a finite
simplicial set.  A finite simplicial set can only have finitely many
non-degenerate simplices, while a simplicial finite set need only
be finite degreewise.

There is an inclusion functor $i: \pro s\F \map \pro \sSet$
from pro-simplicial finite sets to pro-simplicial sets.  We next
define its adjoint.

\begin{defn}
\label{defn:Fin}
If $cX$ is any constant pro-space,
then \mdfn{$\Fin cX$} is the system of all simplicial finite
quotients of $X$.  
If $X$ is an arbitrary pro-space,
then \mdfn{$\Fin X$} is $\lim_s^{\pros} \Fin cX_s$,
where the limit is calculated within the category of
pro-simplicial finite sets.
\end{defn}

\begin{lem}
\label{lem:Fin-i}
The functor $\Fin$ is the left adjoint of the inclusion $i$.
\end{lem}

\begin{proof}
We begin by showing that $\Fin cX$ is a cofiltered system of spaces.
Let $X_1$ and $X_2$ be two simplicial finite quotients of $X$.  We need
to find another simplicial finite quotient $X_3$ of $X$ that refines both
$X_1$ and $X_2$.  
Note that $X_1 \times X_2$ is a simplicial finite set but not necessarily
a quotient of $X$ because the canonical map $X \map X_1 \times X_2$ is
not surjective.  Define $X_3$ to be the image of the map 
$X \map X_1 \times X_2$.  Now $X_3$ is a simplicial finite set because
it is a subobject of $X_1 \times X_2$.  The map $X \map X_3$ is surjective
by construction, so $X_3$ is a simplicial finite quotient of $X$.  
The two maps $X \map X_1$ and $X \map X_2$ factor through $X_3$.
This shows that $\Fin cX$ is a cofiltered system.

Next we will show that $\Fin cX$ has the correct adjoint property.
We want to show that $\Hom_{\pros}(\Fin cX, Y)$ is isomorphic to
$\Hom_{\pros}(cX, Y)$ for every pro-simplicial finite set $Y$.  This
follows from the fact that every map from $X$ into a simplicial finite set 
factors through a simplicial finite quotient of $X$.

Now let $X$ be an arbitrary pro-simplicial set.
The construction of limits in
pro-categories given in Section \ref{subsctn:pro-ctgy} implies 
that 
\[
\Hom_{\pros}(\lim\nolimits_a^{\pros} Z^a, cY) =
\colim_a \Hom_{\pros}(Z^a, cY)
\]
for any cofiltered system $a \mapsto Z^a$ of pro-objects and any constant
pro-object $cY$.  
The desired adjointness 
property for $\Fin X$ now follows formally.
\end{proof}

The adjoint functors $\Fin$ and $i$ connect the categories of pro-simplicial
sets and pro-simplicial finite sets.  Now we have to connect
the categories of pro-simplicial finite sets and simplicial pro-finite sets.
As described explicitly in \cite[\S~3]{I1}, there
are functors $F: \pro s\F \map s \pro \F$ and
$G: s \pro \F \map \pro s\F$ such that $G$ is the left adjoint of $F$
and such that the composition $FG$ is naturally isomorphic to the
identity on $s \pro \F$.  This uses the fact that the category
of simplicial finite sets is small and that the simplicial indexing
category $\Delta^{\op}$ has finite morphism sets.
In other words, the category $s \pro \F$ is a retract of
the category $\pro s\F$.  
As observed above, $F$ and $G$ are
not inverse equivalences of categories because $GF$ is not
naturally isomorphic to the identity functor.

The construction of $G$ is complicated; fortunately
we will not need the details here.
For later reference, we describe the functor $F$.  
Let $X$ be a pro-simplicial finite set.  For each $n \geq 0$, 
$X_n$ is a pro-finite set.  Thus, $[n] \mapsto X_n$ is a 
simplicial pro-finite set, and this is $FX$.

In order to pass between pro-simplicial sets and
simplicial pro-finite sets, we use the compositions
$F \circ \Fin$ and $i \circ G$.  Unfortunately, these functors are
the composition of a left adjoint and a right adjoint.  Thus, they
do not have nice adjointness properties.  This means that we will
not be able to produce a Quillen equivalence \cite[Defn.~8.5.20]{Hi} between
$\pro \sSet$ and $s \pro \F$.  

One might hope that there is a 
$\Z/p$-cohomology model structure
on the intermediate category $\pro s\F$.  Then there would be a 
zig-zag of Quillen equivalences
\[
\xymatrix@1{
\pro \sSet \ar@<2pt>[r] & \pro s\F \ar@<2pt>[l]\ar@<2pt>[r] 
   & s \pro \F. \ar@<2pt>[l]   }
\]
However, since the category of simplicial finite sets is not
a model category, the techniques used in this
paper do not seem to apply.  
Possibly there is another approach altogether.

In the absence of a Quillen equivalence, we have to show directly
that the homotopy categories $\Ho(\pro \sSet)$ and
$\Ho(s \pro \F)$ are equivalent.  

\begin{lem}
\label{lem:cohlgy-R-I}
A map $f$ of pro-simplicial finite sets is a $\Z/p$-cohomology
isomorphism if and only if $Ff$ is a $\Z/p$-cohomology isomorphism
of simplicial pro-finite sets.  
A map $g$ of simplicial pro-finite sets is a $\Z/p$-cohomology
isomorphism if and only if $Gg$ is a $\Z/p$-cohomology isomorphism
of pro-simplicial finite sets.  
\end{lem}

\begin{proof}
Let $X$ be a pro-simplicial finite set.  
The cochain complex $C^{*} X$ used to compute
$H^{*} (X; \Z/p)$ is given by $C^{n} X = \colim_s \Hom((X_s)_n, \Z/p)$.
Using the description of the functor $F$ above, we see that 
this is equal to the cochain complex used to compute 
$H^{*} (F X; \Z/p)$.  
This proves the first claim.

For the second claim, let $Y$ be a simplicial pro-finite set.
We want to show that
$Y$ and $G Y$ have naturally isomorphic $\Z/p$-cohomology.
By the previous paragraph, it suffices to compare $FY$ and $FGY$.
Now $FGY$ is isomorphic to $Y$, so we just need to use
the previous paragraph again.
\end{proof}

We have observed that the functor $GF$ is not well-behaved categorically.
Nevertheless, it does have good cohomological properties.

\begin{cor}
\label{cor:cohlgy-R-I}
The counit natural transformation from the functor
$GF$ to the identity functor on $\pro s\F$
is a natural $\Z/p$-cohomology isomorphism.
\end{cor}

\begin{proof}
If $X$ is any pro-simplicial finite set, both parts of the proof of 
Lemma \ref{lem:cohlgy-R-I} imply that $H^*(X; \Z/p)$ is isomorphic
to $H^*(GFX; \Z/p)$.
\end{proof}

\begin{lem}
\label{lem:Fin-aug}
Let $f: X \map Y$ be a map between pro-simplicial sets.  Then $f$
is a $\Z/p$-cohomology isomorphism if and only if $\Fin(f)$ is
a $\Z/p$-cohomology isomorphism.
\end{lem}

\begin{proof}
We will show that for every pro-simplicial set $X$, the natural map
$X \map \Fin X$ is a $\Z/p$-cohomology isomorphism.
Because of Lemma \ref{lem:lim-cohlgy} and the definition of $\Fin$, 
it suffices to assume that $X$ is a simplicial set.
We must show that the natural map
\[
\colim_Y H^{n} (Y; \Z/p) \map H^{n} (X; \Z/p)
\]
is an isomorphism, 
where $Y$ ranges over all simplicial finite quotients of $X$.
To do this, 
we consider reduced cochain complexes given in degree $n$
by functions into $\Z/p$ from the non-degenerate part $NX_n$ 
of $X$ in degree $n$.

To show that the map of reduced cochain complexes
is surjective, consider an arbitrary cochain $\alpha$, which is 
just a function $NX_n \map \Z/p$.  We need to construct a simplicial
finite quotient $X'$ of $X$ 
and a cochain $\alpha'$ on $X'$ that pulls back to $\alpha$.
Begin by defining an $n$-dimensional simplicial set $Y$ 
whose $(n-1)$-skeleton is
trivial and whose non-degenerate $n$-simplices correspond to the elements
of $\Z/p$.  There is an obvious
map $\sk_n X \map Y$ induced by $\alpha$.  Adjointness gives a map
$X \map \cosk_n Y$.  Since $Y$ is a simplicial finite set, so is
$\cosk_n Y$.  Finally, take $X'$ to be the image in $\cosk_n Y$
of $X$.

To show that the map of reduced cochain complexes
is injective, suppose that $X'$ and $X''$ are
two simplicial finite quotients of $X$, and let $\alpha'$ and $\alpha''$ be
reduced cochains on $X'$ and $X''$ respectively that pull back to
the same reduced cochain $\alpha$ on $X$.  There exists a simplicial finite
quotient $Y$ of $X$ refining both $X'$ and $X''$ (see the proof of
Lemma \ref{lem:Fin-i}).
We now have the diagram
\[
\xymatrix{
X \ar[rrd]\ar[rd]\ar[rdd] \\
& Y \ar[r]\ar[d] & X' \ar[d]^{\alpha'} \\
& X'' \ar[r]_{\alpha''} & \Z/p           }
\]
in which the outer quadrilateral and the two triangles are commutative.
We want to show that the square is also commutative.  This follows 
from the fact that the map $X \map Y$ is surjective.
\end{proof}

\begin{prop}
The functor $F \circ \Fin$ induces a functor
\[
\Ho (\pro \sSet) \map \Ho(s \pro \F)
\]
on homotopy categories, 
and the functor $i \circ G$ induces a functor
\[
\Ho (s \pro \F) \map \Ho(\pro \sSet)
\]
on homotopy categories.
\end{prop}

\begin{proof}
By the universal property of localizations of categories, it suffices
to show that the two functors take weak equivalences to weak equivalences.
Let $f$ be any $\Z/p$-cohomology isomorphism of pro-simplicial sets.  
Lemmas \ref{lem:cohlgy-R-I} and \ref{lem:Fin-aug} imply that
$F \circ \Fin f$ is a $\Z/p$-cohomology isomorphism.
Hence, $F \circ \Fin$ preserves weak equivalences.
For $i \circ G$, this is the second part of Lemma \ref{lem:cohlgy-R-I}.
\end{proof}

\begin{thm}
\label{thm:compare}
The functors $F \circ \Fin$ and $i \circ G$ induce inverse equivalences
between the homotopy categories $\Ho (\pro \sSet)$ and
$\Ho (s \pro \F)$.
\end{thm}

\begin{proof}
The composition $(F \circ \Fin) \circ (i \circ G)$ is isomorphic to
the identity
because $\Fin \circ i$ is the identity by construction of $\Fin$
and because $F \circ G$ is isomorphic to the identity.
On the other hand,
Corollary \ref{cor:cohlgy-R-I} and Lemma \ref{lem:Fin-aug} 
tell us that for every pro-simplicial set
$X$, there are natural weak equivalences
\[
\xymatrix@1{
X \ar[r]^-{\sim} & \Fin X & G F \circ \Fin X \ar[l]_-{\sim}.    }
\]
Thus $X$ and $(i \circ G) \circ (F \circ \Fin) X$ are naturally
isomorphic in $\Ho(\pro \sSet)$.
\end{proof}

\section{Free Groups}
\label{sctn:free-gp}

The point of this section is to describe the both the 
cohomological and homological
$\Z/p$-completions of the 
Eilenberg-Mac Lane space $K(F_n, 1)$, where $F_n$ is the free group on
$n$ generators.
This means that we need to find a fibrant replacement for $cK(F_n,1)$
in the $\Z/p$-cohomological model structure and
the $\Z/p$-homological model structure.  As we will see below,
these two completions turn out to be the same.

Part of the definition of these fibrant replacements requires that the
pro-space be strictly fibrant.
For the rest of this section, we will drop this
requirement.  This change preserves the homotopy type of each space
in the cofiltered system
because of the nature of strict weak equivalences.  Thus, for calculational
purposes we do not really need the strict fibrancy.

\begin{thm}
\label{thm:free-gp}
Consider
the system $\{ K(F_n/H,1) \}$ as $H$ ranges over all normal
subgroups of $F_n$ such that $F_n/H$ is a finite $p$-group; the
structure maps are induced by the canonical quotient maps.  
This pro-space is the fibrant replacement for $K(F_n,1)$ in
either the $\Z/p$-cohomological or $\Z/p$-homological model structure.
\end{thm}

To be precise, we really should also produce a map $cK(F_n,1) \map X$
that induces an isomorphism in $\Z/p$-cohomology.  We will not worry
about this because the map will be obvious and natural in everything
that we do.

\begin{proof}
For notational convenience, write $X_H$ for the space $K(F_n/H,1)$.
We need to show that
$H^*(X;M)$ is
isomorphic to $H^*(F_n;M)$ for every $\Z/p$-module $M$.  
This will show that $X$ and $K(F_n,1)$ are weakly equivalent in
both the $\Z/p$-cohomological and $\Z/p$-homological model structures.
We also have to show that each space $X_H$
is nilpotent in the sense of Definition \ref{defn:K-nilp} 
with respect to the class
of Eilenberg-Mac Lane spaces of the form $K(\Z/p,n)$.
This will show that $X$ is fibrant in both model structures.

First of all, the diagram $X$ 
is a cofiltered system because $F_n/(H \cap K)$ is a finite $p$-group
whenever $F_n/H$ and $F_n/K$ are finite $p$-groups.

Since $K(F_n,1)$ is just
a wedge of $n$ circles, $H^q(F_n; \Z/p)$ is isomorphic to $\Z/p$ in
dimension 0; to $(\Z/p)^n$ in dimension 1; and to the zero group otherwise.
This tells us exactly what the cohomology of $X$ should be.

Let $Z$ be any connected space whose homotopy groups are finite $p$-groups.
Then $Z$ is $\Z/p$-nilpotent in the sense of Bousfield and Kan.  This
can be proved by showing that if $G$ is a finite $p$-group acting on
another finite $p$-group $A$, then $G$ acts nilpotently.
If in addition $Z$ has only finitely many non-zero homotopy groups,
then $Z$ must be $K$-nilpotent in the sense of 
in the sense of Definition \ref{defn:K-nilp}
because of \cite[Prop.~III.5.3]{BK} and Proposition \ref{prop:nilp-tower}.

Since each $X_H$ satisfies the hypotheses in the previous paragraph,
we conclude that each $X_H$ is $K$-nilpotent.
It remains to calculate the $\Z/p$-cohomology of $X$.  This is done
below in Lemma \ref{lem:gp-cohlgy}.
\end{proof}

\begin{lem}
\label{lem:zero-map}
There exists a normal subgroup $H$ of $F_n$
such that $F_n/H$ is a finite $p$-group and such that the map
$H^1(F_n; M) \map H^1(H; M)$ is the zero map for all $\Z/p$-modules $M$.
\end{lem}

\begin{proof}
By an argument similar to the one given after Definition \ref{defn:R-cohlgy},
it suffices to consider the case $M = \Z/p$.
Let $H$ be the kernel of the homomorphism
$F_n \map (\Z/p)^n$ that is the composition of abelianization with
reduction modulo $p$.
Now $K(H,1) \map K(F_n,1)$ is a covering map of degree $p^n$.  More concretely,
$K(H,1)$ is the Cayley graph of the group $(\Z/p)^n$ relative to the
standard basis.  Thus, $K(H,1)$ has one vertex for each element of
$(\Z/p)^n$. The edges of $K(H,1)$ are of the form $v$ to $v + e_i$,
where $v$ is any element of $(\Z/p)^n$ and $e_i$ is any element of
the standard basis.

We use a wedge of $n$ circles as our model for $K(F_n,1)$.
Let $\alpha$ be a $1$-cocycle on $K(F_n,1)$ whose value on
the $i$th circle of $K(F_n,1)$ is the element $\alpha_i$ of $\Z/p$.
Let $\beta$ be the $1$-cocycle on $K(H,1)$ induced by $\alpha$.
The value of $\beta$ on the edge from $v$ to $v+e_i$ is $\alpha_i$.

We construct a $0$-cocycle $\gamma$ whose coboundary is $\beta$.
Let the value of $\gamma$ on the vertex $(v_1, \ldots, v_n)$ of $K(H,1)$
equal $\alpha_1 v_1 + \cdots + \alpha_n v_n$.  Thus, $\beta$
is zero in cohomology, which means that the desired map is zero.
\end{proof}

For each normal subgroup $H$ of $F_n$ such that $F_n/H$ is
a finite $p$-group, we have a short exact sequence
\[
H \map F_n \map F_n/H
\]
which gives rise to a fiber sequence
\[
K(H,1) \map K(F_n,1) \map K(F_n/H,1).
\]
Each such sequence has an associated cohomological Serre spectral sequence
\[
E_2^{st} = H^s(F_n/H; \H^t(H;M)) \Rightarrow
H^{s+t} (F_n; M),
\]
where $\H^t(H; M)$ is a local system on $K(F_n/H,1)$.
Since the Serre spectral sequence is natural
and since filtered colimits respect filtrations, we can take
colimits everywhere and get another spectral sequence
\[
E_2^{st} = \colim_H H^s(F_n/H; \H^t(H;M)) \Rightarrow
H^{s+t} (F_n; M).
\]
Recall how the structure maps of the colimit in the above formula 
are constructed.
If $K$ is a subgroup of $H$, let $\pi$ be the projection
$F_n/K \map F_n/H$.
The map
\[
H^s(F_n/H; \H^t(H;M)) \map H^s(F_n/K; \H^t(K;M)) 
\]
is the composition 
\[
H^s(F_n/H; \H^t(H;M)) \map H^s(F_n/K; \pi^* \H^t(H;M)) \map
H^s(F_n/K; \H^t(K;M)),
\]
where $\pi^* \H^t(H;M)$ is a pullback of local systems.

\begin{lem}
The group
$\colim_H H^s(F_n/H; \H^t(H;M))$ is zero unless $t = 0$.
\end{lem}

\begin{proof}
For $t \geq 2$, each local system $\H^t(H;M)$ is zero because $H$
is a free group.
It only remains to consider the case $t=1$.  
Because $H$ is free,
Lemma \ref{lem:zero-map} 
implies that there exists a subgroup $K$ such that the map
\[
\pi^* \H^1(H;M) \map \H^1(K;M)
\]
of local systems is zero.  This gives the desired result for $t=1$.
\end{proof}

\begin{lem}
\label{lem:gp-cohlgy}
The map $\colim_H H^q(F_n/H;M) \map H^q(F_n;M)$ is
an isomorphism, where the colimit ranges over all normal subgroups of
$F_n$ such that $F_n/H$ is a finite $p$-group.
\end{lem}

\begin{proof}
By the previous lemma, the $E_2$-term of the Serre spectral sequence described
above is concentrated on the line $t = 0$.  This gives the desired
isomorphism.
\end{proof}

\section{Questions}
\label{sctn:?}

The work in this paper leaves some obvious further questions unanswered.
We mention a few of these here in the interest of encouraging future work
on the subject.

\begin{ques}
Are the model structures of 
Theorems \ref{thm:R-cohlgy-mc} and \ref{thm:R-hlgy-mc} are 
right proper?
\end{ques}

The general machinery of localizations does not automatically produce 
right proper model structures.  Presumably the Serre spectral sequence
is the way to approach this problem, but one has to deal with
twisted coefficients.

\begin{ques}
If $X$ is a space considered as a constant pro-space, how do its two
fibrant replacements compare?
\end{ques}

We know that the model structures of Theorems \ref{thm:R-cohlgy-mc}
and \ref{thm:R-hlgy-mc} are distinct.  
However, in Section \ref{sctn:free-gp} we showed that the two
fibrant replacements of $K(F_n,1)$ are the same.  It is easy to imagine
that this would generalize to any space $X$ with some kind of finiteness
hypothesis on the cohomology of $X$.

\begin{ques}
If the ground ring $R$ is $\Z/p$, 
what is the difference between the two fibrant replacements
of a pro-space that is an \'etale topological type?
\end{ques}

Certain kinds of pro-spaces are more relevant to applications than others.
The pro-spaces that arise as \'etale topological types of well-behaved schemes
\cite{AM} \cite{F} are particularly interesting.
Perhaps theorems in algebraic geometry about \'etale cohomology 
with finite coefficients can be used to conclude that the two
fibrant replacements are the same.

\begin{ques}
\label{q:nilp}
Let $R$ be an infinite ring.
If $X$ is a space such that each component is nilpotent in the sense
of Bousfield and Kan and such that the size of $\pi_0 X$ is no bigger
than the size of $R$, can we conclude that $X$ is also
nilpotent in the sense of Section \ref{sctn:nilp}?
\end{ques}

Nilpotence in the sense of Section \ref{sctn:nilp} is definitely 
distinct from Bousfield-Kan nilpotence.
If $X$ is nilpotent in the sense of Section \ref{sctn:nilp}, then
each component of $X$ is nilpotent in the sense of Bousfield and Kan,
and $\pi_0 X$ cannot be bigger than $R$.  The question is whether
the implication works in reverse.

\begin{ques}
Let $R$ be a finite ring.
If $X$ is a space such that each component is nilpotent in the sense
of Bousfield and Kan and such that $\pi_0 X$ is finite,
can we conclude that $X$ is also
nilpotent in the sense of Section \ref{sctn:nilp}?
\end{ques}

This question is just a minor variation on Question \ref{q:nilp}.
Again, we know that nilpotence implies the given conditions.

\end{document}